\documentclass[12pt]{article}
\usepackage{amssymb,amsmath}
\usepackage{cases}
\usepackage{amsfonts}
\usepackage{color,xcolor}
\usepackage[left=2.0cm,right=2.0cm,top=2.0cm,bottom=2.0cm]{geometry}
\usepackage[colorlinks,citecolor=blue,urlcolor=blue]{hyperref}

\newtheorem{theorem}{Theorem}[section]

\newtheorem{lemma}{Lemma}[section]

\newtheorem{remark}{Remark}[section]

\newcommand{\bal}{\begin{align}}
\newcommand{\bbal}{\begin{align*}}
\newcommand{\beq}{\begin{equation}}
\newcommand{\eeq}{\end{equation}}
\newcommand{\bca}{\begin{cases}}
\newcommand{\eca}{\end{cases}}
\def\div{\mathord{{\rm div}}}
\newcommand{\pa}{\partial}
\newcommand{\fr}{\frac}
\newcommand{\na}{\nabla}
\newcommand{\De}{\Delta}

\newcommand{\cd}{\cdot}
\newcommand{\ep}{\varepsilon}
\newcommand{\dd}{\mathrm{d}}

\newcommand{\R}{\mathbb{R}}

\newcommand{\les}{\lesssim}

\begin{document}
\title{Global smooth solutions of the 3D Hall-magnetohydrodynamic equations with large data}

\author{Jinlu Li$^{1}$\footnote{E-mail: lijinlu@gnnu.cn}, Xing Wu$^{2}$\footnote{E-mail: ny2008wx@163.com( Corresponding author)}\\
\small $^1$\it School of Mathematics and Computer Sciences, Gannan Normal University, \\
\small Ganzhou, Jiangxi, 341000, China\\
\small $^2$\it College of Information and Management Science,
Henan Agricultural University,\\
\small Zhengzhou, Henan, 450002, China}

\date{}

\maketitle\noindent{\hrulefill}

{\bf Abstract:} In this paper, we establish the global existence to the three-dimensional incompressible
Hall-MHD equations for a class of large initial data, whose $L^{\infty}$ norms can be arbitrarily large. In addition , we give an example to show that such a large initial value does exist. Our idea is splitting the generalized heat equations  from Hall-MHD system to generate a small quantity for large time $t$.

{\bf Keywords:} Hall-MHD; Global existence; Large initial data.

{\bf MSC (2010):} 35Q35; 76D03; 86A10
\vskip0mm\noindent{\hrulefill}

\section{Introduction}\label{sec1}
This paper studies the Cauchy problem for the following 3D incompressible, resistive Hall-MHD equations
\begin{eqnarray}\label{3D-hmhd}
        \left\{\begin{array}{ll}
          \partial_tu+u\cd\na u+\mu(-\De)^{\alpha} u+\na p=b\cd\na b,& x\in \R^3,t>0,\\
          \partial_tb+u\cd\na b-\nu\Delta b+\nabla\times((\na\times b)\times b)=b\cd\na u,& x\in \R^3,t>0,\\
         \div u=\div b=0,& x\in \R^3,t\geq0,\\
          (u,b)|_{t=0}=(u_0,b_0),& x\in \R^3,\end{array}\right.
        \end{eqnarray}
where $u=(u_1(t,x),u_2(t,x),u_3(t,x))$ and $b=(b_1(t,x),b_2(t,x),b_3(t,x))$ denote the velocity field and magnetic field, respectively, $p\in \R$ is the scalar pressure. $\mu$ is the viscosity and $\nu$ is the magnetic diffusivity. $\Lambda=(-\Delta)^{\fr12}$ is the Zygmund operator and the fractional power operator $\Lambda^{\gamma}$ with $0<\gamma<1$ is defined by Fourier multiplier with symbol $|\xi|^{\gamma}$ (see e.g. \cite{Jacob 2005}), namely,
\begin{eqnarray*}
   \Lambda^{\gamma}u(x)=\mathcal{F}^{-1}|\xi|^{\gamma}\mathcal{F}u(\xi).
\end{eqnarray*}

 Comparing with the standard MHD system, the Hall term $\nabla\times((\na\times b)\times b)$ is included due to Ohm's law, which is believed to be a key issue for understanding magnetic reconnection in the case of large magnetic shear and describes many physical phenomena such as magnetic reconnection in space plasmas \cite{Forbes 1991}, star formation \cite{Balbus 2001,Wardle 2004}, neutron stars \cite{Shalybkov 1997} and also geo-dynamo \cite{Mininni 2003}. For the physical background of the magnetic reconnection and the Hall-MHD, we refer the readers to \cite{Forbes 1991,Lighthill 1960,Simakov 2008} and references therein.
 When $\alpha=1$ and the Hall term is neglected, (\ref{3D-hmhd}) reduces to the standard MHD equation which has been extensively studied, and there are a lot of excellent work, see \cite{Li 2017, Lin 2014, Majda 2001, Sermange 1983} and references therein.

 For the standard Hall-MHD equations, Acheritogaray et al. \cite{Acheritogaray 2011} obtained the global weak solutions in the periodic setting by using the Galerkin approximation. Chae, Degond and Liu  \cite{Chae 2014} proved the global existence of weak solutions as well as the local well-posedness of smooth solutions in the whole space. Meanwhile, they showed that if $\|u_0\|_{H^m}+\|b_0\|_{H^m} (m>\frac{5}{2})$ is small enough, the local smooth solution is global in time. Chae and Lee \cite{Lee 2014} improved  their results under weaker smallness assumptions on the initial data. Laterly, the conditions on the initial data were once more refined by Wan and Zhou \cite{Wan 2015}. Recently, Wan, Zhou \cite{Wan 2019} removed the restriction on $\epsilon$ in \cite{Wan 2015}, established the global existence of strong
solution for standard Hall-MHD equations with the Fujita-Kato type initial data.  Global small solutions for generalized  Hall-MHD equations, we refer the readers to \cite{Ye 2015, Pan 2016, Wu2 2017}. The mathematical studies on \eqref{3D-hmhd} have motivated a relatively large number of research papers concerning the local well-posedness \cite{Chae 2015, Wu1 2018}, regularity criterion \cite{He 2016, Dai 2016, Wan1 2019}, asymptotic behavior \cite{Chae 2013, 1Weng 2016, 2Weng 2016, Wu2 2017}and we can refer the readers to the reference therein.

However, there are few results of global existence for Hall-MHD equations with general initial data without smallness conditions. It is also worth to mention that when $\alpha=1$, $b=0$, the system \eqref{3D-hmhd} is reduced to the Navier-Stokes equations. Lei, Lin and Zhou \cite{Lei 2015} constructed a family of finite energy smooth large solutions to the Navier-Stokes equations with the initial data close to a Beltrami flow. Li, Yang and Yu \cite{Li 2019} established a class of global large solution to the 2D MHD equations with damp terms whose initial energy can be arbitrarily large.  But the Hall term heightens the level of nonlinearity of the standard MHD system from a second-order semilinear to a second-order quasilinear level, significantly making its qualitative analysis more difficult. Until very recently, Zhang\cite{Zhang 2019} obtained  global large smooth solutions in the sense that the initial data can be  arbitrarily large in $H^3(\mathbb{R}^3)$. Li et al. \cite{1Li 2019} constructed a class of large solution with spectrally supported $u_0$ and $b_0=-\nabla\times u_0.$  Different from the constructed  large initial data in \cite{Zhang 2019, 1Li 2019}, for some class of large initial data whose $L^\infty$ norms can be arbitrarily large, by splitting the generalized heat equations  from system \eqref{3D-hmhd} to generate a small quantity, the solution of the system  \eqref{3D-hmhd} with $\mu, \nu>0$  evolve into a global solution.

 Our main result is stated as follows.
\begin{theorem}\label{the1.1} Let $0\leq \alpha\leq 1$ and $U_0$ be a smooth function satisfying ${\rm{div}}U_0=0$, $\na\times U_0=\Lambda U_0$ and
\begin{eqnarray}\label{Equ1.2}
\mathrm{supp} \ \hat{U}_0(\xi)\subset\mathcal{C}\triangleq\Big\{\xi\in \R^3 \big| \ 1-\ep\leq  |\xi|\leq 1+\ep\Big\},\quad 0<\ep<\frac{2-\sqrt{2}}{2}.
\end{eqnarray}
Assume that the initial data fulfills $u_0=v_0+U_0$ and $b_0=c_0+ U_0$, then there exists a sufficiently small positive constant $\delta$, and a universal constant $C$ such that if
\begin{align}\label{condition}
(||v_0||^2_{H^3}+||c_0||^2_{H^3}+\ep||U_0||_{L^2}\big(1+||\hat{U}_0||_{L^1}\big)\Big)
\exp\Big( C(||1+||\hat{U}_0||_{L^1})(||\hat{U}_0||_{L^1}+\ep||U_0||_{L^2})\Big)\leq \delta,
\end{align}
then the system \eqref{3D-hmhd} has a unique global solution.
\end{theorem}

\begin{remark}\label{rem1.1}
Assume that $\hat{a}_0: {\mathbb R}^3\to [0, 1]$ be a radial, non-negative, smooth function which is supported in $\mathcal{C}$ and
$\hat{a}_0\equiv 1$ for $1-\frac12\ep\leq  |\xi|\leq 1+\frac12\ep$.

Notice that
$$a_0(x)=\int_{\R^3}\cos(x\cdot\xi)\hat{a}_0(\xi)\dd \xi \quad\mbox{and}\quad \Lambda^{-1} a_0(x)=\int_{\R^3}\cos(x\cdot\xi)\frac{\hat{a}_0(\xi)}{|\xi|}\dd \xi,$$
then we have $a_0,\Lambda^{-1} a_0\in \R$ and also let $U_0=V_0+\Lambda^{-1} \na\times V_0$ with
\begin{eqnarray*}
&V_0=\ep^{-1}\Big(\log\log\frac1\ep\Big)^\frac12\na\times
\begin{pmatrix}
a_0 \\ 0 \\ 0
\end{pmatrix}
=\ep^{-1}\Big(\log\log\frac1\ep\Big)^\frac12
\begin{pmatrix}
0 \\ \pa_3a_0 \\ -\pa_2a_0
\end{pmatrix}.
\end{eqnarray*}
Here, we can verify that $\mathrm{div} U_0=0$ and $\na \times U_0=\Lambda U_0$.

Moreover, we also have
\begin{eqnarray*}
\hat{U}_0=\ep^{-1}\Big(\log\log\frac1\ep\Big)^\frac12
\begin{pmatrix}
\xi^2_2+\xi^2_3 \\ -\xi_1\xi_2+i\xi_3|\xi| \\ -\xi_1\xi_3-i\xi_2|\xi|
\end{pmatrix}\frac{\hat{a}_0(\xi)}{|\xi|}.
\end{eqnarray*}
Then, direct calculations show that
\begin{align*}
||\hat{U}_0||_{L^1}\approx \Big(\log\log\frac1\ep\Big)^\frac12\quad\mbox{and}\quad||{U}_0||_{L^2}\approx \ep^{-\fr12}\Big(\log\log\frac1\ep\Big)^\frac12
\end{align*}
Thus, the left side of \eqref{condition} becomes
\begin{align*}
C\ep^{\frac{1}{2}}\log\log\frac1\ep
\exp\big(C\log\log\frac1\ep\big).
\end{align*}
Therefore, choosing $\ep$ small enough, we deduce that the system \eqref{3D-hmhd} has a global solution. Notice that $U^1_0=-\ep^{-1}\Big(\log\log\frac1\ep\Big)^\frac12(\pa^2_{2}+\pa^2_{3})\Lambda^{-1} a_0$ and $\hat{U}^1_0\geq 0$, we can deduce that
\bbal
||U^1_0||_{L^\infty}\approx ||\hat{U}^1_0||_{L^1}\gtrsim \Big(\log\log\frac1\ep\Big)^\frac12.
\end{align*}
Let $v_0=c_0=0$, then we have $||u_0||_{L^\infty}=||b_0||_{L^\infty}\gtrsim \Big(\log\log\frac1\ep\Big)^\frac12$.
\end{remark}

{\bf Notations}: Let $\beta=(\beta_1,\beta_2,\beta_3)\in \mathbb{N}^3$ be a multi-index and $D^{\beta}=\pa^{|\beta|}/\pa^{\beta_1}_{x_1}\pa^{\beta_2}_{x_2}\pa^{\beta_3}_{x_3}$ with $|\beta|=\beta_1+\beta_2+\beta_3$. For the sake of simplicity, $a\lesssim b$ means that $a\leq Cb$ for some ``harmless" positive constant $C$ which may vary from line to line. $[A,B]$ stands for the commutator operator $AB-BA$, where $A$ and $B$ are any pair of operators on some Banach space $X$. We also use the notation $||f_1,\cdots,f_n||_{X}\triangleq||f_1||_{X}+\cdots+||f_n||_{X}$.

\section{Reformulation of the System}\label{sec2}
\setcounter{equation}{0}
Setting $U=e^{-\mu t}U_0$ and $B=e^{-\nu t}U_0$, we know that $(U,B)$ solve the following system
\begin{eqnarray}\label{2}
        \left\{\begin{array}{ll}
          \pa_t U+\mu(-\Delta)^\alpha U=\mu((-\Delta)^\alpha-\mathbb{I})U:=F,\\
          \pa_t B-\nu\Delta B=\nu(-\Delta-\mathbb{I})B:=G,\\
          \div U=\div B=0,\\
          (U,B)|_{t=0}=(U_0,U_{0}).\end{array}\right.
\end{eqnarray}

Noticing the fact that $B\cd\na U-U\cd\na B=0$ and denoting $v=u-U$ and $c=b-B$, we can reformulate the system \eqref{3D-hmhd} and \eqref{2} equivalently as
\begin{eqnarray}\label{5}
        \left\{\begin{array}{ll}
\partial_tv+v\cd\na v-c\cd \na c-\De v+\na \Big(p+\frac{|U|^2-|B|^2}{2}\Big)=f+f_1-F,\\
\partial_tc+v\cd\na c-c\cd\na v-\De c+\nabla\times((\na\times c)\times c)=g+g_1+g_2-G,\\
\div v=\div c=0,\\
(v,c)|_{t=0}=(v_0,c_0),\end{array}\right.
\end{eqnarray}
where
\bbal
&f=B\cd\na B-U\cd\na U-\na\Big(\frac{|B|^2-|U|^2}{2}\Big)=(\Lambda B-B)\times B-(\Lambda U-U)\times U,
\\&f_1=B\cd\na c+c\cd\na B-U\cd\na v-v\cd\na U,
\\&g=-\na\times((\na\times B)\times B)=-\na\times((\Lambda B-B)\times B),
\\&g_1=B\cd\na v+c\cd\na U-U\cd\na c-v\cd\na B,
\\&g_2=-\nabla\times((\na\times c)\times B)-\nabla\times((\na\times B)\times c).
\end{align*}

\section{Proof of Theorem \ref{the1.1}}\label{sec3}
\setcounter{equation}{0}
In this section, we present the proof of Theorem \ref{the1.1}. Before proceeding on, we present some estimates which will be used in the proof of Theorem \ref{the1.1}.

\begin{lemma}\label{lem3.1} For small enough $\ep$, under the assumptions of Theorem \ref{the1.1}, the following estimates hold
\bal\label{estimate-l}
||F||_{H^3}+||G||_{H^3}\leq C\max\{\mu,\nu\}e^{-\min\{\mu,\nu\}t}\ep||U_0||_{L^2}
\end{align}
and
\bal\label{estimate-2}
||f||_{H^3}+||g||_{H^3}\leq Ce^{-\min\{\mu,\nu\}t}\ep||U_0||_{L^2}||\hat{U}_0||_{L^1}.
\end{align}
\end{lemma}
{\bf Proof of Lemma \ref{lem3.1}}\quad
For the term $f$, we can show that
\bbal
||F||^2_{H^3}=\mu^2e^{-2\mu t}\int_{1-\ep\leq |\xi|\leq 1+\ep}(1+|\xi|^2)^3\big||\xi|^{2\alpha}-1\big|^2|\hat{U}_0|^2d\xi\leq C\mu^2e^{-2\mu t}\ep^2\alpha^2||U_0||^2_{L^2}.
\end{align*}
Similar argument as the term $F$, we also have
\bbal
||G||^2_{H^3}\leq  C\nu^2e^{-2\nu t}\ep^2||U_0||^2_{L^2}.
\end{align*}

Using the classical Kato-Ponce product estimates and the fact the Fourier transform of a distribute belonging to $L^1$ lies in $L^\infty$, after a simple calculation, we obtain
\bbal
||f||_{H^3}&\lesssim~ ||U||_{L^\infty}||(\Lambda -\mathbb{I})U||_{H^{3}}+||U||_{H^3}||(\Lambda -\mathbb{I})U||_{L^\infty}
\\&\quad + ||B||_{L^\infty}||(\Lambda -\mathbb{I})B||_{H^{3}}+||B||_{H^3}||(\Lambda -\mathbb{I})B||_{L^\infty}
\\&\leq Ce^{-\min\{\mu,\nu\}t}\ep||U_0||_{L^2}||\hat{U}_0||_{L^1},
\end{align*}
and
\bbal
||g||_{H^3}&\lesssim~ ||B||_{L^\infty}||(\Lambda -1)B||_{H^{4}}+||B||_{H^4}||(\Lambda -1)B||_{L^\infty}
\\&\leq Ce^{-\min\{\mu,\nu\}t}\ep||U_0||_{L^2}||\hat{U}_0||_{L^1}.
\end{align*}
Thus, we complete the proof of Lemma \ref{lem3.1}. $\Box$\\

{\bf Proof of Theorem \ref{the1.1}}\quad   Applying $D^\beta$ on $\eqref{5}_1$ and $\eqref{5}_2$ respectively and taking the scalar product of them with $D^\beta v$ and $D^\beta c$ respectively, adding them together and then summing the result over $|\beta|\leq 3$, we get
\bal\label{z0}
\fr12\frac{\dd}{\dd t}\Big(||v||^2_{H^3}+||c||^2_{H^3}\Big)+||\Lambda^\alpha v||^2_{H^3}+||\na c||^2_{H^3}\triangleq\sum^{11}_{i=1}I_i,
\end{align}
where
\bbal
&I_1=-\sum_{0<|\beta|\leq 3}\int_{\R^3}[D^{\beta},v\cd] \na v\cd D^\beta v\dd x
-\sum_{0<|\beta|\leq 3}\int_{\R^3}[D^{\beta},v\cd] \na c\cd D^\beta c\dd x,
\\&I_2=\sum_{0<|\beta|\leq 3}\int_{\R^3}[D^{\beta},c\cd] \na c\cd D^\beta v\dd x
+\sum_{0<|\beta|\leq 3}\int_{\R^3}[D^{\beta},c\cd] \na v\cd D^\beta c\dd x,
\\&I_3=\sum_{0<|\beta|\leq 3}\int_{\R^3}D^{\beta}((\na\times c)\times c)\cd  D^{\beta}(\nabla\times c)\dd x,
\\&I_4=-\sum_{0<|\beta|\leq 3}\int_{\R^3}D^{\beta}(U\cd \na v)\cd D^\beta v\dd x-\sum_{0<|\beta|\leq 3}\int_{\R^3}D^{\beta}(U\cd \na c)\cd D^\beta c\dd x,
\\&I_5=\sum_{0<|\beta|\leq 3}\int_{\R^3}D^{\beta}(B\cd \na c)\cd D^{\beta}v\dd x+\sum_{0<|\beta|\leq 3}\int_{\R^3}D^{\beta}(B\cd \na v)\cd D^{\beta}c\dd x,
\\&I_6=\sum_{0\leq|\beta|\leq 3}\int_{\R^3}D^{\beta}(c\cd \na B)\cd D^{\beta}v\dd x-\sum_{0\leq|\beta|\leq 3}\int_{\R^3}D^{\beta}(v\cd \na B)\cd D^{\beta}c\dd x,
\\&I_7=\sum_{0\leq|\beta|\leq 3}\int_{\R^3}D^{\beta}(c\cd \na U)\cd D^{\beta}c\dd x-\sum_{0\leq|\beta|\leq 3}\int_{\R^3}D^{\beta}(v\cd \na U)\cd D^{\beta}v\dd x,
\\&I_8=\sum_{0<|\beta|\leq 3}\int_{\R^3}D^{\beta}((\na\times c)\times B)\cd D^{\beta}(\nabla\times c)\dd x,
\\&I_9=\sum_{0\leq|\beta|\leq 3}\int_{\R^3}D^{\beta}((\na\times B)\times c)\cd D^{\beta}(\nabla\times c)\dd x,
\\&I_{10}=\sum_{0\leq |\beta|\leq 3}\int_{\R^3}D^{\beta}(f-F)\cd D^{\beta}v\dd x+\sum_{0\leq |\beta|\leq 3}\int_{\R^3}D^{\beta}(g-G)\cd D^{\beta}c\dd x.
\end{align*}
Next, we need to estimate the above terms one by one.

According to the commutate estimate (See \cite{Majda 2001}),
\bal\label{l0}
\sum_{|\alpha|\leq m}||[D^{\alpha},\mathbf{g}]\mathbf{f}||_{L^2}\leq C(||\mathbf{f}||_{H^{m-1}}||\na \mathbf{g}||_{L^\infty}+||\mathbf{f}||_{L^\infty}||\mathbf{g}||_{H^m}),
\end{align}
we obtain
\bal
I_1\leq&~\sum_{0<|\beta|\leq 3}||[D^{\beta},v\cd] \na v||_{L^2}||\na v||_{H^2}+\sum_{0<|\alpha|\leq 3}||[D^{\beta},v\cd]\na c||_{L^2}||\na c||_{H^2}\nonumber\\
\leq&~C||\na v||_{L^\infty}||v||_{H^3}||\na v||_{H^2}+C||v||_{H^3}||\na c||^2_{H^2}\nonumber\\
\leq&~C||v||_{H^3}\Big(||\na v||^2_{H^2}+||\na c||^2_{H^2}\Big)\leq C||v||_{H^3}\Big(||\Lambda^\alpha v||^2_{H^3}+||\na c||^2_{H^3}\Big),\label{z1}\\
I_2\leq&~\sum_{0<|\alpha|\leq 3}||[D^{\alpha},c\cd] \na c||_{L^2}||\na v||_{H^2}+\sum_{0<|\alpha|\leq 3}||[D^{\alpha},c\cd] \na v||_{L^2}||\na c||_{H^2}\nonumber\\
\leq&~C||\na c||_{H^2}||\na v||_{H^2}||c||_{H^3}\nonumber\\
\leq&~C||c||_{H^3}\Big(||\na v||^2_{H^2}+||\na c||^2_{H^2}\Big)\leq C||c||_{H^3}\Big(||\Lambda^\alpha v||^2_{H^3}+||\na c||^2_{H^3}\Big).\label{z2}
\end{align}
Using the cancelation equality $[D^{\beta}(\na\times c)\times c]\cd D^{\beta}(\nabla \times c)=0$ and \eqref{l0}, we get
\bal
I_3=&~-\sum_{0<|\beta|\leq 3}\int_{\R^3}[D^{\beta},c\times](\na\times c)\cd D^{\beta}(\nabla\times c)\dd x\nonumber\\
\leq&~C||\na c||_{L^\infty}||c||_{H^3}||\na c||_{H^3}\nonumber\\
\leq&~ C||c||^2_{H^3}||\na c||^2_{H^3}+\frac18||\na c||^2_{H^3}.\label{z3}
\end{align}
Invoking the following calculus inequality which is just a consequence of Leibniz's formula,
\bbal
\sum_{|\beta|\leq 3}||[D^{\beta},\mathbf{g}]\mathbf{f}||_{L^2}\leq C(||\na \mathbf{g}||_{L^\infty}+||\na^3 \mathbf{g}||_{L^\infty})||\mathbf{f}||_{H^2},
\end{align*}
we obtain
\bal
I_4\leq&~\sum_{0<|\beta|\leq 3}||[D^{\beta},U\cd] \na v||_{L^2}||\na v||_{H^2}+\sum_{0<|\beta|\leq 3}||[D^{\beta},U\cd] \na c||_{L^2}||\na c||_{H^2}\nonumber\\
\leq&~C\Big(||\na U||_{L^\infty}+||\na^3 U||_{L^\infty}\Big)\Big(||v||^2_{H^3}+||c||^2_{H^3}\Big)\leq C||U||_{L^\infty}\Big(||v||^2_{H^3}+||c||^2_{H^3}\Big),\label{z4}\\
I_5\leq&~\sum_{0<|\beta|\leq 3}||[D^{\beta},B\cd] \na c||_{L^2}||\na v||_{H^2}+\sum_{0<|\beta|\leq 3}||[D^{\beta},B\cd] \na v||_{L^2}||\na c||_{H^2}\nonumber\\
\leq&~C\Big(||\na B||_{L^\infty}+||\na^3 B||_{L^\infty}\Big)\Big(||v||^2_{H^3}+||c||^2_{H^3}\Big)\leq C||B||_{L^\infty}\Big(||v||^2_{H^3}+||c||^2_{H^3}\Big),\label{z5}\\
I_8\leq&~\sum_{0<|\beta|\leq 3}||[D^{\beta},B\times] (\na\times c)||_{L^2}||\na c||_{H^3}\nonumber\\
\leq&~C\Big(||\na B||_{L^\infty}+||\na^3 B||_{L^\infty}\Big)||\na c||_{H^2}||\na c||_{H^3}\nonumber\\
\leq&~C\Big(||\na B||_{L^\infty}+||\na^3 B||_{L^\infty}\Big)^2||c||^2_{H^3}+\frac18||\na c||^2_{H^3}\leq C||B||^2_{L^\infty}||c||^2_{H^3}+\frac18||\na c||^2_{H^3}.\label{z6}
\end{align}
By Leibniz's formula and H\"{o}lder's inequality, one has
\bal
I_6\leq&~||c\cd \na B||_{H^3}||v||_{H^3}+||v\cd \na B||_{H^3}||c||_{H^3}\nonumber\\
\leq&~ C\Big(||\na B||_{L^\infty}+||\na^4 B||_{L^\infty}\Big)\Big(||v||^2_{H^3}+||c||^2_{H^3}\Big)\leq  C||B||_{L^\infty}\Big(||v||^2_{H^3}+||c||^2_{H^3}\Big),\label{z7}\\
I_7\leq&~||c\cd \na U||_{H^3}||c||_{H^3}+||v\cd \na U||_{H^3}||v||_{H^3}\nonumber\\
\leq&~ C\Big(||\na U||_{L^\infty}+||\na^4 U||_{L^\infty}\Big)\Big(||v||^2_{H^3}+||c||^2_{H^3}\Big) \leq C||U||_{L^\infty}\Big(||v||^2_{H^3}+||c||^2_{H^3}\Big),\label{z8}\\
I_9\leq&~ ||(\na\times B)\times c||_{H^3}||\na c||_{H^3}\nonumber\\
\leq&~ C\Big(||\na B||_{L^\infty}+||\na^4 B||_{L^\infty}\Big)||c||_{H^3}||\na c||_{H^3}\nonumber\\
\leq&~C\Big(||\na B||_{L^\infty}+||\na^4 B||_{L^\infty}\Big)^2||c||^2_{H^3}+\frac18||\na c||^2_{H^3}\leq C||B||^2_{L^\infty}||c||^2_{H^3}+\frac18||\na c||^2_{H^3}.\label{z9}
\end{align}
Using  H\"{o}lder's inequality and Young inequality, we deduce
\bal
I_{10}&\leq~ C||f,g,F,G||_{H^3}(||v||_{H^3}+||c||_{H^3})\nonumber
\\&\leq~ C||f,g,F,G||_{H^3}+C||f,g,F,G||_{H^3}(||v||^2_{H^3}+||c||^2_{H^3}).\label{z10}
\end{align}
Taking all the estimates \eqref{z1}--\eqref{z10} into \eqref{z0}, we obtain
\bal\label{z12}
&\quad \frac{\dd}{\dd t}\Big(||v||^2_{H^3}+||c||^2_{H^3}\Big)+||\Lambda^\alpha v||^2_{H^3}+||\na c||^2_{H^3}\nonumber\\&\les
\Big(||v||_{H^3}+||c||_{H^3}+||c||^2_{H^3}\Big)\Big(||\Lambda^\alpha v||^2_{H^3}+||\na c||^2_{H^3}\Big)
\nonumber\\&\quad +\Big(||B,U||_{L^\infty}+||B||^2_{L^\infty}\Big)\Big(||v||^2_{H^3}+||c||^2_{H^3}\Big)+||f,g,F,G||^2_{H^3}
\nonumber\\&\les
\Big((||v||^2_{H^3}+||c||^2_{H^3})^{\fr12}+||c||^2_{H^3}\Big)\Big(||\Lambda^\alpha v||^2_{H^3}+||\na c||^2_{H^3}\Big)
\nonumber\\&\quad +\Big(||B,U||_{L^\infty}+||B||^2_{L^\infty}+||f,g,F,G||_{H^3}\Big)\Big(||v||^2_{H^3}+||c||^2_{H^3}\Big)+||f,g,F,G||_{H^3}.
\end{align}
The Hausdorff-Young inequality together with the estimates for $f, g, F, G$ in Lemma \ref{lem3.1} yields
\bal\label{z15}
&\quad \frac{\dd}{\dd t}\Big(||v||^2_{H^3}+||c||^2_{H^3}\Big)+||\Lambda^\alpha v||^2_{H^3}+||\na c||^2_{H^3}\nonumber\\&\leq
C\Big((||v||^2_{H^3}+||c||^2_{H^3})^{\fr12}+||c||^2_{H^3}\Big)\Big(||\Lambda^\alpha v||^2_{H^3}+||\na c||^2_{H^3}\Big)
\nonumber\\&\quad+Ce^{-\min\{\mu,\nu\}t}\Big(1+||\hat{U}_0||_{L^1}\Big)\Big(||\hat{U}_0||_{L^1}+\ep||U_0||_{L^2}\Big)\Big(||v||^2_{H^3}+||c||^2_{H^3}\Big)
\nonumber\\&\quad+Ce^{-\min\{\mu,\nu\}t}\ep||U_0||_{L^2}\Big(1+||\hat{U}_0||_{L^1}\Big).
\end{align}
Now, we define
\bbal
\Gamma\triangleq\sup\{t\in[0,T^*):\sup_{\tau\in[0,t]}\Big(||v(\tau)||^2_{H^3}+||c(\tau)||^2_{H^3}\Big)\leq \eta\},
\end{align*}
where $\eta$ is a small enough positive constant which will be determined later on.

Assume that $\Gamma<T^*$. For all $t\in[0,\Gamma]$, we obtain from \eqref{z15} that
\bbal
\frac{\dd}{\dd t}\Big(||v||^2_{H^3}+||c||^2_{H^3}\Big) &\leq Ce^{-\min\{\mu,\nu\}t}\Big(1+||\hat{U}_0||_{L^1}\Big)\Big(||\hat{U}_0||_{L^1}+\ep||U_0||_{L^2}\Big)\Big(||v||^2_{H^3}+||c||^2_{H^3}\Big)
\nonumber\\&\quad+Ce^{-\min\{\mu,\nu\}t}\ep||U_0||_{L^2}\Big(1+||\hat{U}_0||_{L^1}\Big),
\end{align*}
which follows from the assumption that
\bbal
||v||^2_{H^3}+||c||^2_{H^3}\leq &C\Big(||v_0||^2_{H^3}+||c_0||^2_{H^3}+\ep||U_0||_{L^2}\big(1+||\hat{U}_0||_{L^1}\big)\Big)\\
&\;\exp\Big( C(||1+||\hat{U}_0||_{L^1})(||\hat{U}_0||_{L^1}+\ep||U_0||_{L^2})\Big)\leq C\delta.
\end{align*}
Choosing $\eta=2C\delta$, thus we can get
\bbal
\sup_{\tau\in[0,t]}\Big(||v(\tau)||^2_{H^3}+||c(\tau)||^2_{H^3}\Big)&\leq \fr\eta2 \quad\mbox{for}\quad t\leq \Gamma.
\end{align*}

So if $\Gamma<T^*$, due to the continuity of the solutions, we can obtain there exists $0<\epsilon\ll1$ such that
\bbal
\sup_{\tau\in[0,t]}\Big(||v(\tau)||^2_{H^3}+||c(\tau)||^2_{H^3}\Big)&\leq \fr\eta2 \quad\mbox{for}\quad t\leq \Gamma+\epsilon<T^*,
\end{align*}
which is contradiction with the definition of $\Gamma$.

Thus, we can conclude $\Gamma=T^*$ and
\bbal
\sup_{\tau\in[0,t]}\Big(||v(\tau)||^2_{H^3}+||c(\tau)||^2_{H^3}\Big)&\leq C<\infty \quad\mbox{for all}\quad t\in(0,T^*),
\end{align*}
which implies that $T^*=+\infty$. This completes the proof of Theorem \ref{the1.1}. $\Box$

\section*{Acknowledgments} J. Li is supported by the National Natural Science Foundation of China (Grant No.11801090).

\end{document}